\documentclass[12pt]{amsart}

\usepackage{amssymb,amsmath,amsxtra,amscd, feynmf}
\usepackage{pstricks}
\usepackage{pst-tree, pst-node}
\usepackage[mathscr]{eucal} 
\theoremstyle{plain}
\numberwithin{equation}{section}
\newtheorem{theorem}{Theorem}

\numberwithin{lemma}{section}
\numberwithin{theorem}{section}
\numberwithin{corollary}{section}
\numberwithin{proposition}{section}

\theoremstyle{definition}
\newtheorem{definition}{Definition}
\numberwithin{definition}{section}

\theoremstyle{plain}

\numberwithin{example}{section}
\newtheorem{remark}{Remark}

\newcommand{\nc}{\newcommand}
\nc{\C}{\mathcal{C}}
\nc{\CC}{\widetilde{C}}
\nc{\wt}{\overline}
\nc{\mc}{\mathcal}
\nc{\on}{\operatorname}
\nc{\Cl}{Cl}
\nc{\drva}{\widehat{\OOmega}}
\nc{\spec}{\OOn{Spec}}
\nc{\AutO}{\OOn{Aut} \mc{O}}
\nc{\vac}{|0\rangle}
\nc{\Z}{\mathbb{Z}}
\nc{\zf}[1]{z^{\frac{1}{#1}}}
\nc{\wf}[1]{w^{\frac{1}{#1}}}
\nc{\Mt}{M^{\sigma}}

\nc{\gr}{\OOn{gr}}
\nc{\T}{\mathbb{T}}
\nc{\LT}{\mathbb{LT}}
\nc{\CT}{\mathbb{C}\{ \mathbb{T} \}}
\nc{\g}{\mathfrak{g}}
\nc{\A}{\mathcal{A}}
\nc{\n}{\mathfrak{n}}
\nc{\I}{\mathcal{I}}
\nc{\Hc}{\mathcal{H}}
\nc{\HH}{\mathbf{H}}
\nc{\M}{M}
\nc{\U}{\mathcal{U}}
\nc{\F}{\mathcal{F}}
\nc{\OO}{\mathcal{O}}
\nc{\LRF}{\mathcal{LRF}}
\nc{\FD}{\mathcal{LFG}}
\nc{\LFD}{\mathcal{LFG}}
\nc{\FG}{\mathcal{LFG}}
\nc{\LFG}{\mathcal{LFG}}
\nc{\labb}{\emph{lab}}

\nc{\al}{\alpha}
\nc{\OOl}{\OOverline}

\begin{document}

\title{Feynman graphs, rooted trees, and Ringel-Hall algebras}
\author{Kobi Kremnizer } \thanks{}
%\footnotetext{\kern-15.3pt AMS Subject Classification: 17B10, 17B70, 17B05, 17B81.}
\address{Department of Mathematics, MIT}
\email{kreminze@math.mit.edu}
\author{Matt Szczesny}
\address{Department of Mathematics  
         Boston University, Boston, MA 02215}
\email{szczesny@math.bu.edu}

%\date{April 2008}

\begin{abstract}
We construct symmetric monoidal categories $\LRF, \FD$ of rooted forests and Feynman graphs. These categories closely resemble finitary abelian categories, and in particular, the notion of Ringel-Hall algebra applies. The Ringel-Hall Hopf algebras of $\LRF, \FD$, $\HH_{\LRF}, \HH_{\FD}$ are dual to the corresponding Connes-Kreimer Hopf algebras on rooted trees and Feynman diagrams. We thus obtain an interpretation of the Connes-Kreimer Lie algebras  on rooted trees and Feynman graphs as Ringel-Hall Lie algebras.
%As an application we obtain a description of the enveloping algebra of Feynman graphs analogous to the Grossman-Larson Hopf algebra in the case of rooted trees.   
\end{abstract}

\maketitle 

\section{Introduction}

The Connes-Kreimer Hopf algebras on rooted trees and Feynman graphs  $\mc{H}_{\T}, \mc{H}_{FG}$, introduced in \cite{K}, \cite{CK}, describe the algebraic structure of the BPHZ algorithm in the renormalization of perturbative quantum field theories. If we let $\T$ denote the set of (non-planar) rooted trees, and $\mathbb{Q}\{ \T \}$ the $\mathbb{Q}$--vector space spanned by these, then as an algebra, $\Hc_{\T} = \on{Sym}(\mathbb{Q}\{  \T \})$, and the coalgebra structure is given by the coproduct
\[
\Delta(T) = \sum_{C \textrm{ simple cut }} P_C(T) \otimes R_{C} (T)
\]  
where $P_C(T)$ is the forest of branches resulting from the cut $C$, and $R_C(T)$ is the root component remaining "above" the cut (see \cite{CK} for a more detailed definition). 

$\Hc_{FG}$ is defined analogously, with Feynman graphs in place of rooted trees. More precisely, given a perturbative QFT, and denoting by $\mathbb{Q}\{ \Gamma \}$ the vector space spanned by the one-piece irreducible graphs of the theory ($1\on{PI}$ graphs), $\Hc_{FG} = \on{Sym}(\mathbb{Q}\{ \Gamma \})$ as an algebra. Its coalgebra structure is given by 
\[
\Delta(\Gamma) = \sum_{\gamma \in \Gamma} \gamma \otimes \Gamma / \gamma
\]
where the sum is over all (not necessarily connected) subgraphs of $\Gamma$, and $\Gamma / \gamma$ denotes the graph obtained from $\Gamma$ by shrinking each connected component of $\gamma$ to a point. 

$\Hc_{T}$ and $\Hc_{FG}$ are graded commutative Hopf algebras, and so by the Milnor-Moore theorem, their graded duals $\Hc^{*}_{T}$ and $\Hc^{*}_{FG}$ are the isomorphic to the universal enveloping algebras $\U(\n_T)$, $\U(\n_{FG})$ of the nilpotent Lie algebras $\n_T$, $\n_{FG}$ of their primitive elements. We refer to $\n_T$ and $\n_{FG}$ as the Connes-Kreimer Lie algebras on rooted trees and Feynman graphs respectively. 

In this paper, we present a categorification of the Hopf algebras $\U(\n_T)$, $\U(\n_{FG})$, by showing that they arise naturally as the Ringel-Hall algebras of certain categories $\LRF$, $\FD$ of labeled rooted forests and Feynman graphs respectively. We briefly recall the notion of Ringel-Hall algebra. Given an abelian category $\A$ satisfying the finiteness properties $|\on{Hom}(M,N)| < \infty$ and $|\on{Ext}^1(M,N)| < \infty$ (such an abelian category is called \emph{finitary}), the Ringel-Hall algebra of $\A$, $\HH_{\A}$, is the $\mathbb{Q}$--vector space $$\mathbb{Q}\{ [M] \}$$ spanned by the isomorphism classes $[M] \in \A$. It becomes an associative algebra under the product
\[
[M] \times [N] = \sum \frac{g^{L}_{M,N}}{a_M a_N} [L]
\]
where $g^L_{M,N}$ is the number of short exact sequences
\[
0 \rightarrow M \rightarrow L \rightarrow N \rightarrow 0
\]
and $a_M= | \on{Aut}(M)|$.  $\HH_{\A}$ also possesses a coproduct (see section \ref{Hallalgebras}), and
%\[
%\Delta([M]) = \sum_{N \subset M} [N] \otimes [M/N]
%\]
is in fact a cocommutative Hopf algebra, and so isomorphic to the universal enveloping algebra $\U(\n_{\A})$ of the Lie algebra of its primitive elements $\n_{\A}$, called the \emph{Ringel-Hall Lie algebra of} $\A$. 

The categories $\LRF$ and $\FD$ are not abelian, or even additive, but possess all the necessary properties to define the corresponding Ringel-Hall Hopf algebras. These are enumerated in section \ref{categorydef}. We prove that $\n_T \cong \n_{\LRF}$ and $\n_{FG} \cong \n_{\FD}$. 

This paper is organized as follows. In section \ref{Hallalgebras} we review the notion of Ringel-Hall algebra of a finitary abelian category. Section \ref{trees} introduces some terminology relating to rooted trees and forests, and recalls the Connes-Kreimer Lie algebra on rooted trees. In section \ref{categorydef} we construct the category $\LRF$ of labeled rooted forests, and describe some of its properties. The following section (\ref{HallLRF}) applies the notion of Ringel-Hall algebra to $\LRF$ to obtain $\U(\n_T)$. Finally, in sections \ref{FGcat}, \ref{FGCatt} we construct the category $\FD$ in an analogous manner, and show that its Ringel-Hall algebra is isomorphic to $\U(\n_{FG})$. 

\bigskip
\noindent {\bf Acknowledgements:} M.S. would like to thank Dirk Kreimer and Valerio Toledano-Laredo for valuable conversations and for helpful suggestions. 

\section{Ringel-Hall algebras associated to finitary abelian categories}
\label{Hallalgebras}
In this section, we briefly recall the construction of the Ringel-Hall algebra associated to a finitary abelian category. The notion was introduced in \cite{R}, and our treatment borrows heavily from \cite{S}, where we refer the reader for details and proofs. Recall that a small abelian category $\A$ is called \emph{finitary} if: 
\begin{align}
\textrm{ i) For any two objects } M,N \in \on{Ob}(\A) \textrm{ we have } |Hom(M,N)| < \infty \\
\textrm{ii) For any two objects } M,N \in \on{Ob}(\A) \textrm{ we have } |Ext^1 (M,N)| < \infty
\end{align}
For $M,N,L \in \on{Ob}(\A)$, let $G^L_{M,N}$ denote the set of all exact sequences of the form 
\[
0 \rightarrow M \rightarrow L \rightarrow N \rightarrow 0
\]
By $(i)$, $G^L_{M,N}$ is a finite set. Let $g^L_{M,N} = | G^L_{M,N} |$ and $a_M = |Aut(M)|$. As a $\mathbb{Q}$--vector space 
\[
\HH_{\A} := \bigoplus_{[M] \in I(\A)} \mathbb{Q} \{ [M] \}
\]
where $I(\A)$ denotes the set of isomorphism classes of objects in $\A$, and $[M]$ the isomorphism class of the object $M$. $\HH_{\A}$ is an associative algebra with respect to the product
\begin{equation} \label{Hallproduct}
[M] \times [N] = \sum_{[L] \in I(\A)} \frac{g^{L}_{M,N}}{a_M a_N} [L]
\end{equation}
which is finite by property $(ii)$ of $\A$. This product clearly counts the "number" of extensions of $N$ by $M$ up to isomorphism. A more geometric way of expressing this product is as follows. 
Let $\F(\A)$ denote the vector space of  $\mathbb{Q}$--valued functions on $\A$ supported on finitely many isomorphism classes, i..e
\[
\F(\A) := \{ f: I(\A) \rightarrow \mathbb{Q} | |supp(f)| < \infty \}
\]
and let $\on{Fl}^2 (M)$ denote the space parametrizing flags of length two in $\A$
\[
0 = M_0 \subset M_1 \subset M_2 = M \; \; \; \; \;  M_i \in \A
\]
(note that this is the same as a short exact sequence of objects in $\A$). $\F(\A)$ is equipped with a convolution product: for $f,g \in \F(\A)$, let
\begin{equation} \label{geoHallProduct}
f \times g (M) := \int_{\on{Fl}^2 (M)} f(M_1/M_0) g(M_2/M_1)
\end{equation}
Identifying the symbol $[M]$ with the characteristic functions $\delta_M$ of the isomorphism class of $M \in \A$, we see that the product $\delta_M \times \delta_N$ corresponds to $[M] \times [N]$, and so $\F(\A) = \HH_{\A}$ as associative algebras.
In this formulation, the algebra $\F(\A)$ possesses a natural coproduct
\[
\Delta: \F(\A) \rightarrow \F(\A) \otimes \F(\A)
\]
\begin{equation}
\Delta(f)(M,N) := f(M \oplus N) 
\end{equation}
which endows it with the structure of a cocommutative bialgebra. The primitive elements of $\F(\A)$ are those functions supported on indecomposable elements of $\A$, and form a Lie algebra $\mathfrak{n}_{\A}$. $\F(\A)$ can be naturally identified with the universal enveloping algebra $U(\mathfrak{n}_{\A})$.

%More generally, one may consider the Hall algebra attached to an exact category (see \cite{Hu}). 

\section{Rooted trees and forests} \label{trees}

Let $\T$ denote the set of rooted trees. An element $T \in \T$ is a tree (finite, one-dimensional contractible simplicial complex), with a distinguished vertex $r(T)$, called the \emph{root} of $T$. Let $V(T)$ and $E(T)$ denote the set of vertices and edges of $T$,  and let 
$$
| T | = \#  V(T) 
$$
A \emph{labeling} of a tree $T$ by a set $S$ is a bijection $S \rightarrow V(T)$. In what follows, we will frequently consider trees labeled by subsets $S \subset \mathbb{N}$ of the natural numbers.  For example, 

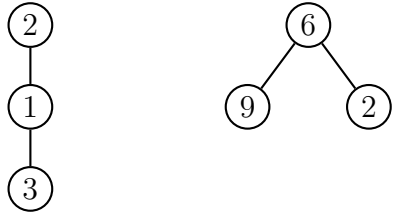
\begin{center} \psset{levelsep=6ex, treesep=1.0cm}
\pstree[]{ \Tcircle{2} }{ \pstree{ \Tcircle{1} } {\Tcircle{3} } } \hspace{2 cm} 
\pstree{\Tcircle{6}}{\Tcircle{9}\Tcircle{2}}
\end{center}
are labeled rooted trees, with the vertex pictured at the top. Let $\LT$ denote the set of rooted trees labeled by subsets of $\mathbb{N}$, and for $T \in \LT$, let $\labb(T) \subset \mathbb{N}$ denote the set of labels (which is canonically identified with $V(T)$). A \emph{labeled rooted forest} $F$ is a set labeled rooted trees, i.e.
\[
F := \{ T_1, T_2, \cdots, T_n \} \hspace{1cm} T_i \in \LT
\]
 An \emph{admissible cut} of a labeled tree $T$ is a subset $C(T) \subset E(T)$ such that at most one member of $C(T)$ is encountered along any path joining a leaf to the root. Removing the edges in an admissible cut divides $T$ into a forest $P_C(T)$ and a rooted tree $R_C(T)$, where the latter is the component containing the root. A \emph{simple cut} is an admissible cut consisting of a single edge. For example, if

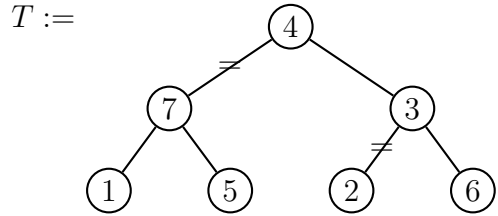
\begin{center} \psset{levelsep=6ex, treesep=1.0cm}
$T:=$ \pstree{  \Tcircle{4}  } {  \pstree{ \Tcircle{7} \ncput{=}   } { \Tcircle{1} \Tcircle{5}}
 \pstree{\Tcircle{3}}{\Tcircle{2}\ncput{=} \Tcircle{6}}   }
\end{center}
and the cut edges are indicated with "=", then 
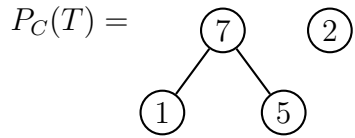
\begin{center}  \psset{levelsep=6ex, treesep=1.0cm}
$P_C(T) = $ \pstree{ \Tcircle{7}}{\Tcircle{1} \Tcircle{5}} \hspace{3cm} \Tcircle{2} 
\end{center}
and
\begin{center} \psset{levelsep=6ex, treesep=1.0cm}
$R_C(T) = $ \pstree{\Tcircle{4}} {\pstree{\Tcircle{3}} {\Tcircle{6}} }
\end{center}
\bigskip
We also allow the \emph{empty} and \emph{full} cuts $C_{null}, C_{full}$, where $$(P_{C_{null}}(T), R_{C_{null}}(T)) = (\emptyset, T) \textrm{ and }  (P_{C_{full}}(T), R_{C_{full}}(T)) = (T, \emptyset)$$ respectively. The latter is considered simple. 
More generally, given a labeled forest $F:=\{ T_1, \cdots, T_n \}$, an admissible cut on $F$ is an $n$--tuple $\{ C_1, C_2, \cdots, C_n \}$ where $C_i$ is an admissible cut of $T_i$, and 
\[
P_C(F) := \{ P_{C_1} (T_1), \cdots, P_{C_n}(T_n) \} \textrm{ and } R_C(F) := \{ R_{C_1}(T_1), \cdots, R_{C_n}(T_n) \} 
\]
\medskip
\begin{definition}
By a \emph{subforest} of a labeled rooted forest $F$ we mean a forest of the form $G =P_C(F)$, for an admissible cut $C$. We write $G \subset F$.
\end{definition}

Given a labeled rooted forest $F$ and two admissible cuts $C_1, C_2$, we write $C_1 < C_2$ if the cut edges of $C_2$ occur closer to the root than those of $C_1$ along any path joining a leaf to the root. Similarly, we write $C_1 \leq C_2$ if the cut edges of $C_2$ occur at those of $C_1$ or closer to the root. The relation $\leq$ defines a partial order on cuts.  We also define
\begin{enumerate}
\item the cut $C_1 \cup C_2$ by the property that $C_i \leq C_1 \cup C_2, \; i=1,2$,  and if $C_i \leq D$ for some cut $D$, then $C_1 \cup C_2 \leq D$.
\item the cut $C_1 \cap C_2$ by the property that $ C_1 \cup C_2 \leq C_i$,  and if $D \leq C_i$ for some cut $D$, then $D \leq C_1 \cap C_2 $.
\end{enumerate}
In other words, $C_1 \cup C_2$ involves cutting the edge closer to the root, and $C_1 \cap C_2$ the farther one. Note that both the operations $\cup, \cap$ are associative and commutative. 
\bigskip

\noindent The following two observations will be important below:
\bigskip
\begin{remark} \label{remark1}
If $G=P_C(F)$ is a subforest of $F$, and $C'$ is an admissible cut on $F$, then $C'$ induces a unique admissible cut on $G$. It is the restriction of the cut $C \cap C'$ to $G$.
\end{remark}
\bigskip
\begin{remark} \label{remark2}
If $G = P_C(F)$ is a subforest of $F$, then there is a bijection between subforests of $R_C(F)$ and subforests $H$ of $F$ such that $G \subset H \subset F$. Both correspond to cuts $C'$ on $F$ such that $C \leq C'$.
\end{remark}
\bigskip
\subsection{The Connes-Kreimer Lie algebra on rooted trees}

In this section, we recall the definition of the Connes-Kreimer Lie algebra on rooted trees $\n_T$ (see \cite{CK}). As a vector space, 
\[
\n_T = \mathbb{Q}\{ \T \} 
\]
i.e. the span of unlabeled rooted trees. On $\n_T$, we have a \emph{pre-Lie} product
"$*$", given, for $T_1, T_2 \in \T$ by
\[
T_1 * T_2 = \sum_{T \in \T} a(T_1,T_2;T) T
\]
where 
\[
a(T_1,T_2;T) := | \{ e \in E(T) | P_{C_e}(T) = T_1, R_{C_e} (T) = T_2 \} |
\]
and $C_e$  denotes the cut severing the edge $e$. 
The Lie bracket on $\n_T$ is given by 
\begin{equation} \label{treebracket}
[T_1,T_2] := T_1 * T_2 - T_2 * T_1
\end{equation}
Thus, for example if
%\begin{align*}
%\psset{levelsep=0.3cm, treesep=0.3cm}
%[\pstree{\Tr{\bullet}}, \pstree{\Tr{\bullet}}{\Tr{\bullet}\Tr{\bullet}} ]&= \pstree{\Tr{\bullet}}{\Tr{\bullet} 
%\Tr{\bullet} \Tr{\bullet }} + \pstree{\Tr{\bullet}}{\pstree{\Tr{\bullet}}{\Tr{\bullet}}\Tr{\bullet}}         - \pstree{\Tr{\bullet}}{\pstree{\Tr{\bullet}}{ \Tr{\bullet} \Tr{\bullet}} } 
%\end{align*}
\begin{center} \psset{levelsep=4ex, treesep=1.0cm}
$T_1:=$ \pstree{\Tcircle{}} {\Tcircle{} {\Tcircle{}} }  \hspace{2cm} \textrm{ and } $T_2:=$ \pstree{\Tcircle{}}{}
\end{center}
then 
\begin{center} \psset{levelsep=4ex, treesep=1.0cm}
$[T_1, T_2] = $ \pstree{\Tcircle{}} {\pstree{\Tcircle{}} {\Tcircle{} \Tcircle{}}} $-$  \pstree{\Tcircle{}}{ {\pstree{\Tcircle{}}{\Tcircle{}}} \Tcircle{}} $-$ \pstree{\Tcircle{}}{\Tcircle{} \Tcircle{} \Tcircle{}}
\end{center}

\section{The category $\LRF$ of labeled rooted forests} \label{categorydef}

Labeled rooted forests can be made into a category $\LRF$ as follows. Let 
\[
\on{Ob}(\LRF) = \{ \textrm{ labeled forests } \} \cup \{ \emptyset \}
\]
where $\emptyset$ denotes the \emph{empty forest}, which plays the role of zero object. 
%By a \emph{full subtree} of a labeled rooted tree $\Lambda \subset T$ we mean a subtree of the form $\Lambda = P_C(T)$, where $C$ is a \emph{simple} cut of $T$.  Note that $\emptyset, T$ are full subtrees of $T$, and $T$ the only full subtree containing the root. 
%By a \emph{full subforest} of a labeled rooted forest $F$ we mean a subforest
%Analogously, we define a \emph{ full subforest} of $F:= \{ T_1, \cdots, T_n \}$ to be a forest $\{ \Lambda_1, \cdots, \Lambda_n \}$, where each $\Lambda_i \subset T_i$ is a full subtree of $T_i$. Equivalently, and more succinctly, full subforests are subforests of the form $P_C (F)$, where $C$ is an admissible cut on $F$.

\begin{definition}
We say that two labeled rooted trees $T_1$ and $T_2$ are \emph{isomorphic}, and write $f:T_1 \cong T_2$ if there exists a root and incidence-preserving bijection $ f: \labb(T_1) \rightarrow \labb(T_2)$. Two forests $F_1 := \{ T_1, \cdots, T_n \}$ and $F_2:= \{ U_1, \cdots, U_n  \}$ are isomorphic if $ f: T_i \cong U_i$ for $i=1, \cdots, n$.  We write $f: F_1 \cong F_2 $
\end{definition}
\bigskip
\noindent If $F_1, F_2 \in \LRF$, we now define 
\begin{align*}
\on{Hom}(F_1,F_2) := & \{ (C_1,C_2,f) | C_i \textrm{ is an admissible cut of } F_i, \\ & \; f: R_{C_1} (F_1) \cong P_{C_2} (F_2) \}
\end{align*}
and the \emph{image}  of $(C_1,C_2,f)$,  $Im(C_1,C_2,f)$ (or $Im(f)$ if $C_1, C_2$ are understood) to be the subtree $P_{C_2}(F_2)$ of $F_2$. For $F \in \LRF$, the morphism $(C_{null}, C_{full}, id)$ is the identity morphism.
%\noindent {\bf Note:} Every morphism determines an admissible cut $C'$ on $F_2$, such that $P_{C'}(F_2)$ is the full subforest isomorphic to $R_C(F_1)$ via $f$. Alternatively, we could have defined morphisms as consisting of the data of two admissible cuts $C, C'$ on $F_1$ and $F_2$ respectively, and an isomorhism $f : R_C(F_1) \rightarrow P_{C'}(F_2)$. 

\bigskip

\noindent {\bf Example:} if 
\bigskip
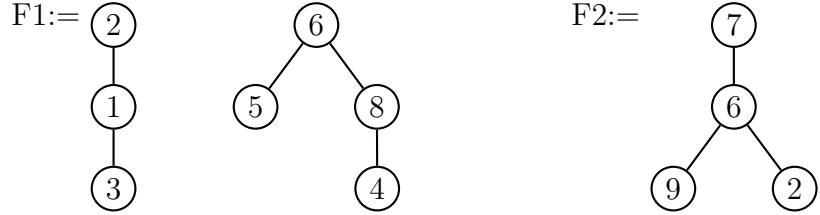
\begin{center} \psset{levelsep=6ex, treesep=1.0cm}
F1:= \pstree[]{ \Tcircle{2} }{ \pstree{ \Tcircle{1} } {\Tcircle{3} } } \hspace{1 cm} 
\pstree{\Tcircle{6}}{\Tcircle{5} {\pstree{\Tcircle{8}}{\Tcircle{4}} }} \hspace{2cm}
F2:= \pstree{\Tcircle{7}}{\pstree{\Tcircle{6}}{\Tcircle{9}\Tcircle{2}}}
\end{center}
\bigskip
then a morphism is given by the triple $(C_1, C_2,f)$ where:
\begin{itemize}
\item $C_1$ is the full cut on the first tree in $F1$, and on the second severs the edge joining $8$ to $4$ 
\item $C_2$ is the cut on $F_2$ which severs the edge joining $7$ to $6$.
\item $f : R_{C_1}(F_1) \cong P_{C_2}(F_2)$ is defined by $f(5) = 9$, $f(6)=6$ ,$f(8)=2$. 
\end{itemize}

\begin{remark}
Note that $(C_1, C_2, f')$, where $f'(5)=2, f'(6)=6, f'(8)=9$ is also a morphism. 
\end{remark}

\bigskip

\noindent The composition of morphisms
\[
\on{Hom}(F_1, F_2) \times \on{Hom}(F_2, F_3) \rightarrow \on{Hom}(F_1,F_3)
\]
is defined as follows. Suppose that $(C_1,C_2,f) \in \on{Hom}(F_1,F_2)$, and $(D_2,D_3,g) \in \on{Hom}(F_2,F_3)$. By remark \ref{remark1}, the cut $D_2$ on $F_2$ induces a cut on the subforest $P_{C_2}(F_2) \cong R_{C_1}(F_1)$, which by remark \ref{remark2} corresponds to a subforest $P_{E_1}(F_1)$ of $F_1$ containing $P_{C_1}(F_1)$. The image $g \circ f (R_{E_1}(F_1)) \subset F_3$ is a subforest $P_{E_3}(F_3)$. We define the composition $(D_2,D_3,g) \circ (C_1,C_2,f)$ to be $(E_1,E_3, g\circ f)$. It is easy to see that this composition is associative. We thus obtain:

\begin{theorem}
With the above definitions of $\on{Ob}(\LRF)$ and $\on{Hom}$, $\LRF$ forms a category.
\end{theorem}

\bigskip

\noindent $\LRF$ has among other, the following properties:
\bigskip
\begin{enumerate}
\item Given labeled forests $F_1, F_2 $ we denote their disjoint union by $F_1 \oplus F_2$. The disjoint union of forests equips $\LRF$ with a symmetric monoidal structure. 
%For $F_1, \cdots, F_n \in \on{Ob}(\LRF)$, $F_1 \oplus \cdots \oplus F_n$ is both a product and a coproduct in $\LRF$. 
\bigskip
\item The empty forest $\{ \emptyset \}$ is an intial, terminal, and null object. 
\bigskip
\item  \label{kernel} Every morphism 
\begin{equation} \label{morphism}
(C_1, C_2,f) : F_1 \rightarrow F_2
\end{equation}
possesses a kernel
\[
(C_{null},C_1, id):  P_C(F_1) \rightarrow F_1 
\]
where $C_{null}$ denotes the empty cut, and $id$ the identity map $id:P_{C_1}(F_1) = R_{C_{null}}(P_{C_1} (F_1)) \cong  P_{C_1}(F_1)$.
\bigskip
\item \label{cokernel} Similarly, every morphism \ref{morphism} possesses a cokernel
\[
(C_2, C_{full}, id) : F_2 \rightarrow R_{C_2} (F_2)
\]
where $id$ is the identity map $R_{C_2}(F_2) \cong R_{C_2}(F_2) = P_{C_{full}}( R_{C_2}(F_2))$.
 
\medskip
\noindent We will frequently use the notation $F_2/F_1$ for $coker((C_1,C_2,f))$. 

\bigskip
\noindent {\bf Note:} Properties \ref{kernel} and \ref{cokernel} imply that the notion of exact sequence makes sense in $\LRF$.
\bigskip
\item All monomorphisms are of the form
\[
(C_{null},C_1, f) : P_{C_1} (F_1) \rightarrow F_1
\]
where $f$ is an automorphism of $P_{C_1}(F_1)$. Once the image subforest $P_{C_1}(F_1)$ is fixed, all monomorphisms with that image form a torsor over $\on{Aut}(P_{C_1})$, and there are therefore $|\on{Aut}P_{C_1}(F_1)|$ of them. 
All epimorphisms are of the form
\[
(C_2,C_{full}, g) : F_2 \rightarrow R_{C_2} (F_2)
\]
where $g$ is an automorphism of $R_{C_2}(F_2)$.The epimorphisms with fixed kernel subtforest $P_{C_2}(F_2)$ form a torsor over $\on{Aut}(R_{C_2}(F_2))$, and so there are $|\on{Aut}(R_{C_2}(F_2))|$ of them. 
\bigskip
\item
Sequences of the form \label{propses}
\begin{equation} \label{ses}
\emptyset \rightarrow P_C(F) \overset{(C_{null},C, id)}{\longrightarrow} F \overset {(C,C_{full},id)}{\longrightarrow} R_C  (F) \rightarrow \emptyset
\end{equation} 
are exact, and it follows from the last property that all other short exact sequences arise by composing with automorphisms of $P_C(F)$ and $R_C(F)$ on the left and right respectively. 
\bigskip
\item \label{quotients} By remark \ref{remark2}, given a forest $F$ and an admissible cut $C$, there is a bijection between subobjects $F'$ of $F$ containing $P_C(F)$, i.e. chains $P_{C}(F) \subset F' \subset F$, and subobjects of $R_C(F)$. 
\bigskip
\item $\on{Hom}(F_1,F_2)$ and $\on{Ext}^{n} (F_1, F_2)$ are finite sets. 
\bigskip
\item We may define the Grothendieck group of $\LRF$, $K(\LRF)$, as 
\[
K(\LRF) = \bigoplus_{[M] \in \on{I}(\LRF)} \mathbb{Z}[M] / \sim
\]
where the equivalence relation $\sim$ is defined by:  $[M] \sim [N]$ iff there exists a short exact sequence 
\[
0 \rightarrow M \rightarrow L \rightarrow N \rightarrow 0
\]
It is then easy to see that $K(\LRF) = \mathbb{Z}$, as every forest is a direct sum of trees, and each tree is an extension of the one-vertex tree. 
\end{enumerate}

\section{The Ringel-Hall algebra of $\LRF$}. \label{HallLRF}

We proceed to define the Ringel-Hall algebra of the category $\LRF$ as in the case of finitary abelian categories. Let $I(\LRF)$ denote the isomorphism classes of objects in $\LRF$, and let 
$$ \HH_{\LRF} := \{ f : I(\LRF) \rightarrow \mathbb{Q} | |supp(f)| < \infty  \}. $$
i.e. the space of $\mathbb{Q}$--valued functions supported on finitely many isomorphism classes in $\LRF$. We equip $\HH_{\LRF} $ with the convolution product
\begin{equation} \label{mult}
f \times g (F) = \sum_{G \subset F} f(G) g(F/G)
\end{equation}
where the notation $F/G$ is used as explained in property \ref{cokernel} of $\LRF$. It is clear that this sum is finite, as any object in $\LRF$ possesses finitely many subobjects. 

The proof of the following theorem is essentially identical to that in the case of finitary abelian categories in \cite{S}. We include it for the sake of completeness. 

\begin{theorem} \label{assocthm}
The multiplication $\times$ in \ref{mult} is associative.
\end{theorem}

\begin{proof}
Suppose $f,g,h \in \HH_{\LRF} $, and $F \in I(\LRF)$. We have
\begin{align*}
(f \times ( g \times h)) (F) &= \sum_{G \subset F} f(G) (g \times h) (F/G) \\
                                    &= \sum_{G \subset F, H' \subset F/G} f(G) g(H') h((F/G)/H') \\
                                    &= \sum_{G \subset H \subset F} f(G) g(H/G) h(F/H)
\end{align*}
\begin{align*}
((f \times g) \times h) (F) &= \sum_{K \subset F} (f \times g)(K) h(F/K) \\
                                  &= \sum_{K \subset F, L \subset K} f(L) g(K/L) h(F/K) \\
                                  &= \sum _{L \subset K \subset F } f(L) g(K/L) h(F/K)
\end{align*}
where the equality between the second and third lines follows from property \ref{quotients} which yields a bijection between the sets
\[
\{ H' \subset F/G  \} \hspace{1.5cm} \leftrightarrow \hspace{1.5cm} \{ G \subset H \subset F \}
\]
satisfying 
\[
H' = H/G \textrm{ and } (F/G)/H' = F/H
\]
\end{proof}
\bigskip
\begin{remark} \label{NecCond}
The only properties of the category $\LRF$ that needed to establish the associativity of the product \ref{mult} are \ref{kernel}, \ref{cokernel},  \ref{propses}
\end{remark}
\bigskip
Property \ref{propses} of $\LRF$ implies that if $g^K_{F_1,F_2}$ is the number of short exact sequences of the form
\begin{equation} \label{ses2}
\emptyset \rightarrow F_1 \rightarrow K \rightarrow F_2 \rightarrow \emptyset
\end{equation}
and $h^K_{F_1,F_2}$ is the number of subobjects $L \subset K$ such that $L \cong F_1$ and $K/L \cong F_2$, then 
\[
h^K_{F_1,F_2} = \frac{g^{K}_{F_1,F_2}}{|\on{Aut}(F_1)||\on{Aut}(F_2)|}
\]
Thus, if $\delta_{F_1}, \delta_{F_2} \in \HH_{\LRF}$ are the characteristic functions of the isomorphism classes of $F_1,F_2$, we have
\[
\delta_{F_1} \times \delta_{F_2} = \sum_{K \in I(\LRF)} h^{K}_{F_1,F_2} \delta_{K} = \sum_{K \in I(\LRF)} \frac{g^{K}_{F_1,F_2}}{|\on{Aut}(F_1)||\on{Aut}(F_2)|} \delta_{K}
\]
The algebra $\HH_{\LRF}$ is graded by $K(\LRF) = \mathbb{Z}$, which coincides with the grading by the number of vertices in a forest. 
We introduce a coproduct on $\HH_{\LRF}$, as in the case of a finitary abelian category, by
\[
\Delta : \HH_{\LRF} \rightarrow \HH_{\LRF} \otimes \HH_{\LRF}
\]
\begin{equation} \label{LRFcoproduct}
\Delta(f)(F,G ) = f(F \oplus G)
\end{equation}

\begin{theorem} \label{HisEnv}
$\HH_{\LRF}$ is a co-commutative Hopf algebra isomorphic to $\U(\n_T)$.
\end{theorem}
\begin{proof}
$\Delta$ is co-commutative since $F \oplus G = G \oplus F$, and the uniqueness of decompositions into indecomposable objects (labeled trees) in $\LRF$ yields coassociativity. It is very easy to check that $\Delta$ is compatible with $\times$. $\HH_{\LRF}$ is therefore a graded connected bialgebra, and thus a Hopf algebra.  The Milnor-Moore theorem implies that $\HH_{\LRF}$ is isomorphic to the universal enveloping algebra of its primitive elements $\n_{\LRF}$, which for the coproduct $\Delta$  are exactly the indecomposable elements of $\LRF$ - the characteristic functions supported on the the isomorphism class of a single labeled tree. 

It remains to verify that for $\delta_{T_1}, \delta_{ T_2 } \in \n_{\LRF}$, the bracket 
$$[\delta_{T_1}, \delta_{T_2}]_{\times} := \delta_{T_1} \times \delta_{T_2 } - \delta_{ T_2 } \times \delta_{ T_1 }$$
coincides with the Lie bracket \ref{treebracket} under the map 
\[
j: \T \rightarrow \HH_{\LRF}
\] 
\[
j(T) = \delta_{T}
\]
(this makes sense since any two labelings of a tree are isomorphic). Extending $j$ linearly, it is easy to see that for unlabeled rooted trees $T_1,T_2$,
\[
j( T_1 \star T_2) = \delta_{T_1} \times \delta_{T_2} - \delta_{T_1 \oplus T_2}
\]
which implies that
\[
j( T_1 \star T_2 - T_2 \star T_1) =  \delta_{T_1} \times \delta_{T_2} -  \delta_{T_2} \times \delta_{T_1}
\]
This proves the result.
\end{proof}

\begin{remark}
The Hopf algebra $\HH_{\LRF}$ is canonically isomorphic to the Grossman-Larson Hopf algebra (see \cite{GL} - the fact that the Grossman-Larson Hopf algebra is isomorphic to $U(\n_T)$ was first proved in \cite{P}, with certain inaccuracies corrected in \cite{H}). For instance, for the forest $F = \{ T_1, \cdots, T_n \} \in I(\LRF)$, 
\[
\Delta(\delta_{F}) = \sum_{J \subset [n] = \{ 1, \cdots, n \}} \delta_{F_J} \otimes \delta_{F_{[n] \backslash J}}
\]
where if $J=\{ j_1, \cdots , j_k \} \subset [n]$, $F_J := \{ T_{j_1}, \cdots, T_{j_k} \} \in I(\LRF)$. The Grossman-Larson product involves a summation over all subtree attachments which easily seen to correspond to the enumeration of all exact sequences \ref{ses2}.

\end{remark}

\section{Feynman graphs} \label{FGcat}

In this section we show how to equip Feynman graphs with the structure of a category $\FD$ possessing properties completely parallel to those of $\LRF$, in such a way that 
$$\HH_{\FD} \cong \U(\n_{FG}),$$ where $\n_{FG}$ is the Connes-Kreimer Lie algebra on Feynman graphs. We thus arrive at an interpretation of the latter as the Ringel-Hall Lie algebra of $\FD$.  Our treatment of the combinatorics of graphs is taken from \cite{Y}. In order to not get bogged down in notation, we focus on the special case of $\phi^3$ theory (the case of trivalent graphs with only one edge-type). The results of this section extend to the general case in a completely straighforward manner. We begin with a series of definitions.

\begin{definition}
A \emph{graph} $\Gamma$ consists of a set $H = H(\Gamma)$ of half-edges, a set $V = V(\Gamma)$ of vertices, a set of vertex-half edge adjacency relations $( \subset V \times H)$, and a set of half edge - half edge adjacency relations $(\subset H \times H)$, with the requirements that each half edge is adjacent to at most one other half edge and to exactly one vertex. Note that graphs may not be connected. 

Half edges which are not adjacent to another half edge are called \emph{external edges}, and denoted $Ex=Ex(\Gamma) \subset E = E(\Gamma)$. Pairs of adjacent half edges are called \emph{internal edges}, and denoted $Int(\Gamma)$. 
\end{definition}

\begin{definition}
A \emph{half edge $S$--labeled graph}, (\emph{labeled graph} for short), is a triple $(\Gamma, S, \rho)$, where $\Gamma$ is a graph, $S$ is a set such that $|S|=|H|$, and  $\rho: H \rightarrow  S $ is a bijection. $S$ will usually be obvious from context. 
\end{definition}

\begin{definition}
\begin{enumerate}
\item A \emph{Feynman graph} is a graph where each vertex is incident to exactly three half-edges, and each connected component has $2$ or $3$ external edges.  We denote the set of Feynman graphs by $FG$. 

\item Similarly, we can define  a \emph{labeled Feynman graph}. We denote the set of labeled Feynman graphs by $LFG$.

\item A graph (or a labeled graph) is \emph{1-particle irreducible} ($1\on{PI}$) if it is connected, and remains connected under the removal of an arbitrary internal edge. 
\end{enumerate}
\end{definition}

\noindent {\bf Example:} The graph $\Gamma_{eg}$

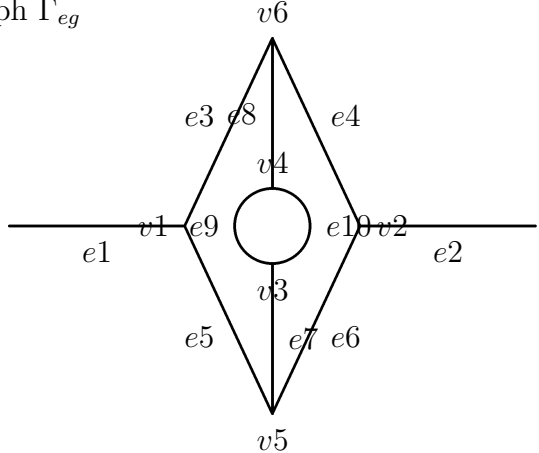
\begin{center} \label{graph1}
\unitlength=1mm
\begin{fmffile}{fig1} 
\begin{fmfgraph*}(70,50)
\fmfleft{i}
\fmfright{o}
\fmftop{t}
\fmfbottom{b}
\fmflabel{$v5$}{b}
\fmflabel{$v6$}{t}
\fmflabel{$v1$}{v1}
\fmflabel{$v2$}{v2}
\fmflabel{$v3$}{v3}
\fmflabel{$v4$}{v4}
\fmf{plain, label=$e1$}{i,v1}
\fmf{plain, label=$e2$}{v2,o}
\fmf{plain, label=$e3$}{t,v1}
\fmf{plain, label=$e4$}{t,v2}
\fmf{plain, label=$e5$}{b,v1}
\fmf{plain, label=$e6$}{b,v2}
\fmf{plain, label=$e7$}{b,v3}
\fmf{plain, label=$e8$}{t,v4}
% \fmf{plain,left,tension=1}{v3,v4,v3}
\fmf{plain, left, label=$e9$, tension=1}{v3,v4}
\fmf{plain, right, label=$e10$, tension=1}{v3,v4}
\end{fmfgraph*}
\end{fmffile}
\end{center}
\vspace{1cm}

Is a $1\on{PI}$ Feynman graph with two external edges. We have labeled each vertex and edge, and each half-edge can be thought of as labeled by a pair $(v,e)$ where $v$ is a vertex, and $e$ is an edge incident to $v$. 

\begin{definition}
Given a Feynman graph $\Gamma$, a \emph{subgraph} $\gamma$ is a Feynman graph such that  $V(\gamma) \subset V(\Gamma)$, $H(\gamma) \subset H(\Gamma)$, and such that if $v \in V(\gamma)$, and $(v,e) \in V(\Gamma) \times H(\Gamma)$, then $e \in H(\gamma)$ (i.e. the subgraph has to contain all half-edges incident to its vertices). We also insist that $\on{dim}_{\mathbb{Q}} (H_{1}(\gamma, \mathbb{Q})) > 0$ (i.e. that a subgraph contain at least one loop). We write $\gamma \subset \Gamma$.  The same definition applies to labeled graphs. 
\end{definition}
\bigskip
\noindent {\bf Example:} We define a subgraph $\gamma_{eg} \subset \Gamma_{eg}$ as follows. Let $V(\gamma_{eg}) = \{ v3, v4 \}$, $E(\gamma_{eg}) := \{ \textrm{ all half-edges incident to } v3, v4 \}$, i.e. $$E(\gamma_{eg}) = \{ (v4,e8), (v3,e7), (v3,e9), (v3,e10), (v4,e9), (v4,e10)  \}$$ and all incidences inherited from $\Gamma_{eg}$. 
\bigskip

\noindent We proceed to define the contraction of subgraphs of Feynman graphs. 

\begin{definition}
Let $\Gamma$ be a Feynman graph, and $\gamma \subset \Gamma$ a connected subgraph. The \emph{quotient graph} $\Gamma / \gamma$ is defined as follows. If $\gamma$ has $3$ external edges, then $\Gamma / \gamma$ is the Feynman graph with 
\begin{enumerate}
\item $V(\Gamma/ \gamma)$ set the vertex set of $\Gamma$ with all vertices of $\gamma$ removed, and a new trivalent vertex $v$ added. 
\item $H(\Gamma/ \gamma)$ the half edge set of $\Gamma$, with all half edges corresponding to internal half edges of $\gamma$ removed.
\item All adjacencies inherited from $\Gamma$, and the external half edges of $\gamma$ joined to $v$.  
\end{enumerate}
\bigskip
If $\gamma$ has $2$ external edges, then $\Gamma/ \gamma$ is the Feynman graph with
\bigskip
\begin{enumerate}
\item $V(\Gamma / \gamma)$ is $V(\Gamma)$ with all the vertices of $\gamma$ removed.
\item $H(\Gamma / \gamma)$ is $H(\Gamma)$ with all half edges of $\gamma$ removed. 
\item All adjacencies inherited from $\Gamma$, as well as the adjacency of the external half-edges of $\gamma$. 
\end{enumerate}
\bigskip
Finally, If $\gamma \subset \Gamma$ is an arbitrary (not necessarily connected) Feynman subgraph, then $\Gamma / \gamma$ is defined to be the Feynman graph obtained by performing successive quotients by each connected component. Note that the order of collapsing does not matter. 
\end{definition}
\bigskip

\noindent {\bf Example:} With $\Gamma_{eg}, \gamma_{eg} $ as above, $\Gamma_{eg} / \gamma_{eg}$ is:
\vspace{1cm}
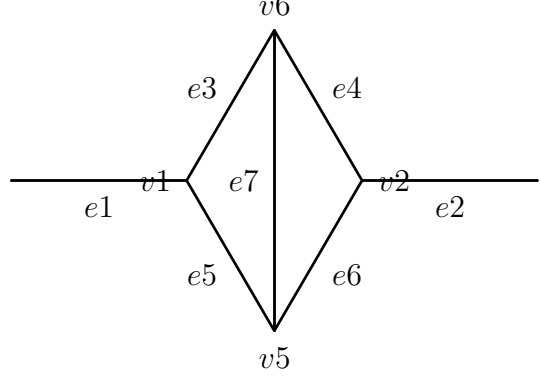
\begin{center} \label{fig2}
\unitlength=1mm
\begin{fmffile}{fig2} 
\begin{fmfgraph*}(70,40)
\fmfleft{i}
\fmfright{o}
\fmftop{t}
\fmfbottom{b}
\fmflabel{$v5$}{b}
\fmflabel{$v6$}{t}
\fmflabel{$v1$}{v1}
\fmflabel{$v2$}{v2}
\fmf{plain, label=$e1$}{i,v1}
\fmf{plain, label=$e2$}{v2,o}
\fmf{plain, label=$e3$}{t,v1}
\fmf{plain, label=$e4$}{t,v2}
\fmf{plain, label=$e5$}{b,v1}
\fmf{plain, label=$e6$}{b,v2}
\fmf{plain, label=$e7$}{t,b}
\end{fmfgraph*}
\end{fmffile}
\end{center}
\vspace{1cm}

\begin{remark} \label{remark3}
If $\gamma_1, \gamma_2 \subset \Gamma$ are subgraphs of a (labeled) Feynman graph, then $\gamma_1 \cap \gamma_2$ is a Feynman subgraph of $\gamma_i$ (and of course also $\Gamma$).
\end{remark}
\bigskip
\begin{remark} \label{remark4}
If $\gamma \subset \Gamma$ is a subgraph of a (labeled) Feynman graph, then there is a bijection between subgraphs of $\Gamma/ \gamma$, and subgraphs $\gamma'$ of $\Gamma$ such that $\gamma \subset \gamma' \subset \Gamma$.
\end{remark}
\bigskip

\subsection{The Connes-Kreimer Lie algebra on Feynman graphs}

In order to define the Connes-Kreimer Lie algebra structure on Feynman graphs, we must first introduce the notion of inserting a graph into another graph. Let
$$\mathbb{Q} \{ LFG \} $$ denote the vector space spanned by labeled Feynman graphs. 

\begin{definition}
Let $\Gamma_1, \Gamma_2 \in LFG$. If $\Gamma_1$ has three external edges, $v \in V(\Gamma_2) $, and $f: Ex(\Gamma_1) \rightarrow H(v)$ is a bijection (where $H(v)$ are the labeled half-edges incident to the vertex $v$), then let $\Gamma_2 \circ_{v,f} \Gamma_1$ be the labeled Feynman graph such that
\bigskip
\begin{itemize}
\item $V(\Gamma_2 \circ_{v,f} \Gamma_1) = V(\Gamma_2) \cup V(\Gamma_1) \backslash v$.
\item $H(\Gamma_2 \circ_{v,f} \Gamma_1) = H(\Gamma_1) \cup_f H(\Gamma_2)$ - i.e. the unions of the half-edges of each graph, with the identifications induced by $f$.
\item The adjacencies induced from those of $\Gamma_1$ and $\Gamma_2$. 
\end{itemize}
\bigskip
If $\Gamma_1$ has two external edges, $\{ e_1,e_2 \} \in Int(\Gamma_2) \subset H \times H$, and $f$ is a bijection between $Ex(\Gamma_1)$ and $\{ e_1, e_2 \}$ (there are two of these), then $\Gamma_2 \circ_{e,f} \Gamma_1$ is the labeled Feynman graph such that
\bigskip 
\begin{itemize}
\item $V(\Gamma_2 \circ_{e,f} \Gamma_1) = V(\Gamma_1) \cup V(\Gamma_2)$.
\item $H(\Gamma_2 \circ_{e,f} \Gamma_1) = H(\Gamma_1) \cup H(\Gamma_2)$.
\item The adjacency induced by $f$ as well as those induced from $\Gamma_1$ and $\Gamma_2$.
\end{itemize}
\bigskip

Let $\n_{LFG}$ denote the $\mathbb{Q}$--vector space spanned by \emph{unlabeled} Feynman graphs. Given a labeled Feynman graph $\Gamma$, denote by $\wt{\Gamma}$ the corresponding unlabeled Feynman graph. Thus, 
\[
\n_{LFG} = \mathbb{Q} \{ LFG \} / \sim
\]
where $\Gamma \sim \Gamma' $ iff $\wt{\Gamma} = \wt{\Gamma'}$.
\noindent We now equip $\n_{LFG}$ with the pre-Lie product "$\star$", defined by
\[
\wt{\Gamma_1} \star \wt{\Gamma_2} := \sum_{v \in V(\Gamma_2), \\
 f: Ex(\Gamma_1) \rightarrow H(v)} \wt{\Gamma_2 \circ_{v,f} \Gamma_1 }
\]
if $\Gamma_1$ has three external edges, and
\[
\wt{\Gamma_1} \star \wt{\Gamma_2} := \sum_{e \in Int(\Gamma_2), \\ f: Ex(\Gamma_1) \rightarrow \{ e_1, e_2 \} } \wt{\Gamma_2 \circ_{v,f} \Gamma_1}
\]
if $\Gamma_1$ has two external edges, and extended linearly ( in the above formulas, we first choose an arbitrary labeling of the Feynman graphs).  
 Finally, we can define the Lie bracket on $\n_{LFG}$ by
\begin{equation} \label{FGLiebracket}
[\wt{\Gamma_1},\wt{\Gamma_2}] := \wt{\Gamma_1} \star \wt{ \Gamma_2 } - \wt{\Gamma_2} \star \wt{ \Gamma_1}
\end{equation}

\end{definition}  
\bigskip
\begin{fmffile}{fig6}
\unitlength=1mm
\noindent {\bf Example:} Suppose
\begin{align}
\Gamma_1 &=
\parbox{20mm}{
\begin{fmfgraph}(20,20)
\fmfleft{i1}
\fmfright{o1}
\fmf{plain}{i1,v1}
\fmf{plain}{v2,o1}
\fmf{plain,left,tension=0.7}{v1,v2,v1}
\end{fmfgraph}} &
\Gamma_2 &= 
\parbox{20mm}{
\begin{fmfgraph}(40,40)
\fmfleft{i}
\fmftop{t}
\fmfbottom{b}
\fmf{plain}{i,v1}
\fmf{plain}{t,v2}
\fmf{plain}{b,v3}
\fmf{plain}{v1,v2,v3,v1}
\end{fmfgraph}}
\end{align}
\end{fmffile}

\noindent then 

\unitlength=1mm
\begin{fmffile}{fig7}
\begin{align*}
[\Gamma_1,\Gamma_2] &= 6 \hspace{0.5cm}
\parbox{20mm}{
\begin{fmfgraph}(40,40)
\fmfleft{i}
\fmftop{t}
\fmfbottom{b}
\fmf{plain}{i,v1}
\fmf{plain}{t,v2}
\fmf{plain}{b,v3}
\fmf{plain}{v2,v3}
\fmf{plain}{v3,v1}
\fmf{plain}{v1,v4}
\fmf{plain}{v2,v5}
\fmf{plain,left,tension=0.8}{v4,v5,v4}
\end{fmfgraph}}
& - 12 \hspace{0.5cm} &
\parbox{20mm}{\begin{fmfgraph}(40,40)
\fmfleft{i}
\fmfright{o}
\fmftop{t}
\fmfbottom{b}
\fmf{plain}{i,v1}
\fmf{plain}{v2,o}
\fmf{phantom}{t,v3}
\fmf{phantom}{b,v4}
\fmf{plain}{v1,v3,v2,v4,v1}
\fmf{plain}{v3,v4}
\end{fmfgraph}}
\end{align*}
\end{fmffile}

\noindent $\n_{LFG}$ has an alternative presentation as follows (see \cite{CK2}). For (unlabeled) Feynman graphs $\wt{\Gamma_1}, \wt{\Gamma_2}$, let 
$$
a(\wt{\Gamma_1}, \wt{\Gamma_2}, \wt{\Gamma}) := | \{ \textrm{ subgraphs } \gamma \subset \Gamma | \gamma \cong \Gamma_1,  \Gamma / \gamma = \Gamma_2 \} |
$$
$\n_{LFG}$ can now be equipped with the pre-Lie structure
\begin{equation} \label{secondPLstructure}
\wt{\Gamma_1} \# \wt{\Gamma_2} := \sum_{\wt{\Gamma}} a(\wt{\Gamma_1}, \wt{\Gamma_2}, \wt{\Gamma}) \wt{\Gamma}
\end{equation}
where the sum is taken over all unlabeled Feynman graphs. We thus obtain "another" Lie bracket on $\n_{LFG}$:
\begin{equation} \label{bracket2}
[\wt{\Gamma_1} , \wt{\Gamma_2}] := \wt{\Gamma_1} \# \wt{\Gamma_2} - \wt{\Gamma_2} \# \wt{\Gamma_1}
\end{equation}
The two structures \ref{FGLiebracket} and \ref{bracket2} are shown in \cite{CK2} to be isomorphic via the map
\[
\wt{\Gamma} \mapsto | \on{Aut} (\wt{\Gamma}) | \wt{\Gamma}
\]

\section{The category $\LFG$ of labeled Feynman graphs} \label{FGCatt}

Labeled Feynman graphs form a category $\LFG$ as follows. Let 
\[
\on{Ob}(\LFG) = \{ \textrm{ labeled Feynman graphs } \} \cup \{ \emptyset \}
\]
where $\emptyset$ denotes the \emph{empty Feynman graph}, which plays the role of zero object. Note that objects of $\LFG$ may have several connected components. 

\begin{definition}
We say that two labeled Feynman graphs $\Gamma_1$ and $\Gamma_2$ are \emph{isomorphic} if there exist bijections $f_V: V(\Gamma_1) \rightarrow V(\Gamma_2)$, $f_H: H(\Gamma_1) \rightarrow H(\Gamma_2)$ which induce bijections on all incidences. We write $f: \Gamma_1 \cong \Gamma_2$.  
\end{definition}
\bigskip
\noindent If $\Gamma_1, \Gamma_2 \in \LFG$, we now define
\begin{align*}
\on{Hom}(\Gamma_1, \Gamma_2) := & \{ (\gamma_1, \gamma_2, f) | \gamma_i \textrm{ is a subgraph of } \Gamma_i, \\ & f: \Gamma_1 / \gamma_1 \cong \gamma_2 \}
\end{align*}
For $\Gamma \in \LFG$, $(\emptyset,\Gamma,id)$ is the identity map in $\on{Hom}(\Gamma,\Gamma)$.

\noindent The composition of morphisms in $\LFG$
\[
\on{Hom}(\Gamma_1, \Gamma_2) \times \on{Hom}(\Gamma_2, \Gamma_3) \rightarrow \on{Hom}(\Gamma_1,\Gamma_3)
\]
is defined as follows. Suppose that $(\gamma_1,\gamma_2,f) \in \on{Hom}(\Gamma_1,\Gamma_2)$, and $(\tau_2,\tau_3,g) \in \on{Hom}(\Gamma_2,\Gamma_3)$. By remark \ref{remark3}, the subgraph $\tau_2$ on $\Gamma_2$ induces a subgraph $\tau_2 \cap \gamma_2$ of $\gamma_2 \subset \Gamma_2$, which by remark \ref{remark4} corresponds to a subgraph of $\xi$ of $\Gamma_1$ containing $\gamma_1$. The image $g \circ f (\xi) \subset \Gamma_3$ is a subgraph $\rho \subset \tau_3$. We define the composition $(\tau_2,\tau_3,g) \circ (\gamma_1,\gamma_2,f)$ to be $(\xi,\rho, g\circ f)$. It is easy to see that this composition is associative. We thus obtain:

\begin{theorem}
With the above definitions of $\on{Ob}(\LFG)$ and $\on{Hom}$, $\LFG$ forms a category.
\end{theorem}

\bigskip

\noindent $\LFG$ has a list of properties completely analogous to $\LRF$:
\bigskip
\begin{enumerate}
\item Given labeled Feynman graphs $\Gamma_1, \Gamma_2 $ we denote their disjoint union by $Gamma_1 \oplus \Gamma_2$. This operation equips $\LFG$ with a symmetric monoidal structure. 
%For $\Gamma_1, \cdots, \Gamma_n \in \on{Ob}(\LFG)$, $\Gamma_1 \oplus \cdots \oplus \Gamma_n$ is both a product and a coproduct in $\LFG$. 
\bigskip
\item The empty graph $\{ \emptyset \}$ is an intial, terminal, and null object. 
\bigskip
\item  \label{kernelFG} Every morphism 
\begin{equation} \label{morphismFG}
(\gamma_1, \gamma_2,f) : \Gamma_1 \rightarrow \Gamma_2
\end{equation}
possesses a kernel
\[
(\emptyset,\gamma_1, id):  \gamma_1 \rightarrow \Gamma_1 
\]
where $id$ the identity map $id: \gamma_1 / \emptyset \cong \gamma_1$.
\bigskip
\item \label{cokernelFG} Similarly, every morphism \ref{morphismFG} possesses a cokernel
\[
(\gamma_2, \Gamma_2 / \gamma_2, id) : \Gamma_2 \rightarrow \Gamma_2 / \gamma_2
\]
where $id$ is the identity map $id: \Gamma_2/ \gamma_2 \rightarrow  \Gamma_2/ \gamma_2 $ 
\medskip
\noindent We will frequently use the notation $\Gamma_2/\Gamma_1$ for $coker((\gamma_1,\gamma_2,f))$. 

\bigskip
\noindent {\bf Note:} Properties \ref{kernelFG} and \ref{cokernelFG} imply that the notion of exact sequence makes sense in $\LFG$.
\bigskip
\item All monomorphisms are of the form
\[
(\emptyset,\gamma_1, f) : \gamma_1 \rightarrow \gamma_1
\]
where $f$ is an automorphism of $\gamma_1 \subset \Gamma_1$.
and all epimorphisms are of the form
\[
(\gamma_2,\Gamma_2/ \gamma_2, g) : \Gamma_2 \rightarrow \Gamma_2 / \gamma_2
\]
where $g$ is an automorphism of $\Gamma_2 / \gamma_2$.
\bigskip
\item
Sequences of the form \label{propsesFG}
\begin{equation} \label{sesFG}
\emptyset \rightarrow \gamma \overset{(\emptyset,\gamma, id)}{\longrightarrow} \Gamma \overset {(\gamma,\Gamma / \gamma,id)}{\longrightarrow} \Gamma/\gamma \rightarrow \emptyset
\end{equation} 
are exact, and all other short exact sequences arise by composing with automorphisms of $\gamma$ and $\Gamma/ \gamma$ on the left and right respectively. 
\bigskip
\item \label{quotients} By remark \ref{remark4}, given a labeled Feynman graph $\Gamma$ and $\gamma \subset \Gamma$, there is a bijection between subobjects $\gamma'$ of $\Gamma$ containing $\gamma$, i.e. chains $\gamma \subset \gamma' \subset \Gamma$, and subobjects of $\Gamma/ \gamma$. 
\bigskip
\item $\on{Hom}(\Gamma_1,\Gamma_2)$ and $\on{Ext}^{n} (\Gamma_1, \Gamma_2)$ are finite sets. 
\bigskip
\item We may define the Grothendieck group of $\LFG$, $K(\LFG)$, as \label{KLFG}
\[
K(\LFG) = \bigoplus_{[M] \in \on{I}(\LFG)} \mathbb{Z}[M] / \sim
\]
where the equivalence relation $\sim$ is defined by:  $[M] \sim [N]$ iff there exists a short exact sequence 
\[
0 \rightarrow M \rightarrow L \rightarrow N \rightarrow 0
\]
$K(\LFG) = \mathbb{Z}[\mathcal{P}]$, where $\mathcal{P}$ is the set of isomorphism classes of \emph{primitive Feynman graphs}, which are those Feynman diagrams not containing any proper subgraphs. This follows since every Feynman graph is a union of connected ones, and each connected component is obtained by repeated insertions of primitive graphs. 
\end{enumerate}

We may now proceed exactly as in section \ref{HallLRF} to define the Ringel-Hall algebra of $\LFG$. We define:
\[
\HH_{\LFG} := \{  f: I(\LFG) \rightarrow \mathbb{Q} | |supp(f)| < \infty  \}
\]
We may equip $\HH_{\LFG}$ with the product \ref{mult} and coproduct \ref{LRFcoproduct}. Since the category $\LFG$ satisfies the conditions in \ref{NecCond}, the proof of associativity goes through as in theorem \ref{assocthm}. Finally, the argument of theorem \ref{HisEnv} establishes the following result:

\begin{theorem}
$\HH_{\LFG}$ is a Hopf algebra isomorphic to $U(\n_{LFG})$.
\end{theorem}

\section{Further directions}

Once Lie algebras of Feynman graphs are thought of categorically as Ringel-Hall Lie algebras, interesting and natural questions emerge. In particular, combining the results of this paper with those of \cite{BTL}  yields an interesting link between perturbative quantum fields theory and the geometry of irregular connections on $\mathbb{P}^1$ and their Stokes phenomena. Very briefly, in \cite{BTL}, the authors study isomonodromic families of irregular connections  with values in the Ringel-Hall Lie algebra of a finitary Abelian category $\A$, $\n_{\A}$. The connections are parametrized by families of \emph{stability conditions} on the category $\A$. 

Recall that a stability condition on an abelian category $\A$ is a homomorphism of abelian groups
\[
Z: K(\A) \mapsto \mathbb{C}
\]
(where $K(\A)$ is the Grothendieck group of $\A$ ), such that 
\[
Z(K_{> 0} (\A)) \subset \mathbb{H}
\]
where  $K_{> 0} (\A)$ is the positive cone generated by the classes of nonzero objects, and $\mathbb{H}$ is the upper half plane. Thus, a stability condition is a choice of $Z(M) \in \mathbb{H}$ for each non-zero object $M \in \mathbb{H}$ such that $Z$ is additive across exact sequences. Given $Z$, each object $M$ has a well-defined \emph{phase} 
\[
\phi(M) = \frac{1}{\pi} \on{arg} (Z) \in (0,1)
\]
and $M$ is  said to be \emph{semistable} with respect to $Z$ if every non-zero subobject $N \subset M$ satisfies $\phi(N) \leq \phi(M)$. As $Z$ changes, and passes through "walls", the collection of semistable objects changes. 

Stability conditions certainly make sense for categories $\LRF, \LFG$, even though they are not abelian. For $\LFG$, whose Grothendieck group is generated by primitive Feynman graphs (see property \ref{KLFG} of $\LFG$), a stability condition $Z$ assigns a phase to each primitive graph. For a given $Z$, the semistable graphs will be allowed to contain certain primitive subgraphs, and not others (additively). Changing $Z$ amounts to changing the allowed subgraphs. 
These, and other connections with the results in \cite{BTL} are explored in \cite{KS2}.

\end{document}